\documentclass[10pt,amssymb]{article}

\hoffset \voffset
\usepackage{amsfonts,amssymb,amsmath,latexsym,verbatim,amscd}
\usepackage{amsthm,amsmath,amssymb,amscd}
\usepackage{amssymb}
\usepackage{graphics}
\usepackage{fancyhdr}

\def\à{\`a}
\def\é{\'e}
\def\è{\`e}
\def\ç{\'o}
\def\ù{\`u}
\def\ò{\`o}
\def\ì{\`\i\ }

\newcommand{\qu}{\,\,}

\newcommand{\dvol}[1]{{\rm \,d}vol_{#1}}

\newcommand{\nor}[2]{\|{#1}\|_{#2}}
\newcommand{\al}{\alpha}
\newcommand{\be}{\beta}
\newcommand{\bigo}[1]{\mathcal{O} \left( #1 \right)}

\newcommand{\zi}{z^{i}}

\newcommand{\ep}{\varepsilon}
\newcommand{\ex}{{\rm{e}}}
\newcommand{\te}{\theta}
\newcommand{\gtil}{\tilde{g}}
\newcommand{\A}[1]{A^{(#1)}}
\newcommand{\gep}{g_{\ep}}
\newcommand{\geps}{g^{\ep}}
\newcommand{\uep}{u_{\ep}}
\newcommand{\tip}{{\rm {th}}}
\newcommand{\cip}{{\rm {ch}}}
\newcommand{\sip}{{\rm {sh}}}
\newcommand{\lep}{\mathcal{L}_{\ep}}

\newcommand{\vin}{v_{\infty}}
\newcommand{\epj}{\ep_{j}}
\newcommand{\fhi}{\varphi}
\newcommand{\dd}[2]{\frac{\partial #1}{\partial #2}}

\newtheorem{teor}{Theorem}[section]

\newtheorem{cor}[teor]{Corollary}
\newtheorem{lemma}[teor]{Lemma}
\newtheorem{prop}[teor]{Proposition}

\title{Generalized connected sum construction for constant scalar curvature metrics}

\author{L. Mazzieri~\thanks{Scuola Normale Superiore di
Pisa and Laboratoire d'Analyse et de Math\'ematiques Appliqu\'ees,
Universit\'e Paris 12. E-mail : l.mazzieri@sns.it}}

\begin{document}

\maketitle

\begin{abstract}
In this paper we construct constant scalar curvature metrics on the
generalized connected sum $M = M_1 \, \sharp_K \, M_2$ of two
compact Riemannian manifolds $(M_1,g_1)$ and $(M_2,g_2)$ along a
common Riemannian submanifold $(K,g_K)$, in the case where the
codimension of $K$ is $\geq 3$ and the manifolds $M_1$ and $M_2$
carry the same nonzero constant scalar curvature $S$. In particular
the structure of the metrics we build is investigated and described.
\end{abstract}

\section{Introduction and statement of the result}

Connected sum of solutions of nonlinear problems has revealed to be
a very powerful tool in understanding solutions of many geometric
problems (minimal and constant mean curvature surfaces \cite{Maz-Pac
2}, \cite{Maz-Pac-Pol}, constant scalar curvature metrics
\cite{Joy}, \cite{Maz-Pol-Uhl}, \cite{Maz-Pac 1}, and recently even
Einstein metrics \cite{And}). However, generalized connected sums
along a submanifold have not been addressed so much, probably
because these constructions are less flexible.

\medskip

In this paper we consider the problem of constructing solutions to
the Yamabe equation (i.e. conformal constant scalar curvature
metrics) on the generalized connected sum $M = M_1 \, \sharp_K \,
M_2$ of two compact Riemannian manifolds $(M_1,g_1)$ and $(M_2,g_2)$
along a common (isometrically embedded) submanifold $(K,g_K)$ of
codimension $\geq 3$. We are able to perform this generalized
connected sum under the assumptions that the two initial Riemannian
metrics have the same constant scalar curvature $S$ and the
linearized Yamabe operator about the metrics $g_i$ (i.e. the
operators $\Delta_{g_i} \qu + \qu S/(n-1) $) have trivial kernels,
for $i=1,2$.

\medskip

To put this result in perspective, let us recall the classical
result of Schoen-Yau \cite{Sch-Yau} and Gromov-Lawson \cite{Gro-Law}
which ensures that if the manifolds $M_1$ and $M_2$ carry positive
scalar curvature metrics, then so does the generalized connected sum
$M = M_1 \, \sharp_K \, M_2$ along a submanifold $K$ of codimension
$\geq 3$ and, thanks to the resolution of the Yamabe problem by T.
Aubin and R. Schoen, $M$ can be endowed with a constant positive
scalar curvature metric. This result however does not give the
precise structure of the constant scalar curvature metric one
obtains on the generalized connected sum $M$. In particular, one
would like to know how does the constant scalar curvature metric on
the connected sum looks like in terms of the constant scalar
curvature metric on the summands. Our result does not cover all
cases covered by the above mentioned result but, as it is typical
for most of the gluing results, we have a very precise description
of the metric on the connected sum in terms of the metric on the
summands. Indeed, away from the region where the generalized
connected sum takes place, we obtain metrics on $M$ which are
conformal to the metrics $g_i$ with some conformal factor as close
to the constant function $1$ as we want.

\medskip

In the case of connected sum at points a result analogous to ours
had been obtained by D. Joyce \cite{Joy}. Our strategy is roughly
speaking the same~: we first write down a one dimensional family of
approximate solutions metrics $(g_{\ep})_{\ep \in (0,1)}$ (where the
parameter $\ep$ represent the size of the tubular neighborhood we
excise from each manifold in order to perform the generalized
connected sum), then, we study the linearized scalar curvature
operator about the metric $g_\ep$ and, for all sufficiently small
$\ep$, we find suitable conformal factors $u_{\ep}$ such that the
metrics $\gtil_{\ep} = u_\ep^{\frac{4}{n-2}} \gep$ have constant
scalar curvature $S$ using a simple fixed point argument. Let us now
describe our result more precisely.

\medskip

Let $(M_1, g_1)$ and $(M_2,g_2)$ be two $m$-dimensional compact
Riemannian manifolds with constant scalar curvature $S$, and suppose
that there exists a $k$-dimensional Riemannian manifold $(K,g_K)$
which is isometrically embedded in each $(M_i,g_i)$, for $i=1,2$, $m
\geq 3$, $m-k \geq 3$. We also assume that the normal bundles of $K$
in $(M_i,g_i)$ can be diffeomorphically identified.
Finally, we assume that on both manifolds, the operator
\begin{eqnarray*}
 L_{g_i} : = \Delta_{g_i} + \frac{S}{n-1}
\end{eqnarray*}
is injective.

\medskip

Let $M = M_1 \,\sharp_{K}\, M_2$ be the generalized connected sum of
$(M_1,g_1)$ and $(M_2, g_2)$ along $K$ which is obtained by removing
an $\ep$-tubular neighborhood of $K$ from each $M_i$ and
identifying the two boundaries. \\

Our main result reads~:
\begin{teor}
Under the above assumptions, it is possible to endow $M$ with a
family of constant scalar curvature metrics $\tilde{g}_{\ep}$, $\ep
\in (0, \ep_0)$ whose scalar curvature $S_{\tilde{g}_{\ep}}$ is
constant equal to $S$. In addition, the following holds

\medskip

(i) - The metric $\tilde{g}_{\ep}$ is conformal to the metrics $g_i$
away from a fixed (small) tubular neighborhood of $K$ in $M_i$,
$i=1,2$ for a conformal factor $u_\ep$ which can be chosen so that
\[
\nor{u_{\ep} - 1}{L^{\infty}(M)} \leq c \, \ep^{\frac{n-2}{2} -
\delta},
\]
where $\max\{0, (n-4)/ 2 \} < \delta < (n-2) / 2 $, $n=m-k$ and
$c>0$ does not depend on $\ep$.

\medskip

(ii) - As $\ep$ tends to $0$, the metrics $\tilde{g}_{\ep}$ converge
to
 $g_i$ on compacts of $M_i \setminus \iota_i(K)$, $i=1,2$.
\end{teor}

A typical case where our result applies is when both $(M_1, g_1) =
(M_2,g_2)$ and $K$ is any submanifold of codimension $\geq 3$,
provided the operator $L_{g_i}$ has no nontrivial kernel.

\medskip

There are some main technical differences between our construction
and D. Joyce's construction in the connected sum case. Our
construction seems to be less flexible in the sense that more
hypothesis are needed on the summands to obtain the result. In
particular (so far) the construction only holds when $(K, g_K)$ is
isometrically embedded in both $(M_i,g_i)$ and if this is not the
case it seems harder to construct a reasonable approximate solution
$g_\ep$ to our problem. The second difference comes from the
analysis of the operator $L_{g_\ep}$, the linearized scalar
curvature operator about the metric $g_\ep$. As in the connected sum
case, the derivation of the estimates of the solution of $L_{g_\ep}
\, u = f$ follows from application of the maximum principle.
However, in the generalized connected sum case, the estimates for
the partial derivatives of the solution $u$ are not as nicely
behaved as in the connected sum case. Hopefully, the scalar
curvature equation is a semilinear elliptic equation and hence, the
nonlinear part of this equation only involves the function $u$ and
not its partial derivatives.

\medskip

It is possible to extend our result to the case where $S=0$ relaxing
the fact that the scalar curvature one obtains on the summand is
equal to $0$. Indeed, in this case, the scalar curvature obtained on
$M$ might not be equal to $0$ but will be a constant close to $0$.

\section{Building the metrics}

Let $(K,g_K)$ be a $k$-dimensional Riemannian manifold isometrically
embedded in both the $n$-dimensional Riemannian manifolds
$(M_1,g_1)$ and $(M_2,g_2)$,
\[
\iota_i : K \hookrightarrow M_i
\]
We assume that the isometric map $\iota_1^{-1}\circ\iota_2 :
\iota_1(K) \rightarrow \iota_2(K) $ extends to a diffeomorphism
between the normal bundles of $\iota_i (K)$ in $(M_i,g_i)$, $i=1,2$.
We further assume that the metrics $g_1$ and $g_2$ have the same
constant scalar curvature $S$. In this section our aim is to perform
a generalized connected sum of $(M_1,g_1)$ and $(M_2,g_2)$ along
$(K,g_K)$ and to construct on the new manifold $M = M_1
\,\sharp_{K}\, M_2$ a family of metrics
$(g_{\ep})_{\ep \in (0,1)}$, whose scalar curvature is close to $S$.\\

For a fixed $\ep \in (0,1)$, we describe the generalized connected
sum construction and the definition of the metric $g_\ep$ in local
coordinates, the fact that this construction yields a globally
defined metric will follow at
once.\\

Let $U^k$ be an open set of $\mathbb{R}^k $, $B^{m-k}$ the
$(m-k)$-dimensional open ball ($m-k \geq 3$). For $i=1,2$, $F_i :
U^k \times B^{m-k} \rightarrow W_i \subset M_i$ given by
\[ F_i(z,x) : = \exp^{M_i}_z(x)
\]
defines local Fermi coordinates near the coordinate patches
$F_i(\cdot,0)\left( U \right)\subset \iota_i(K) \subset M_i$. In
these coordinates, the metric $g_i$ can be decomposed as
\begin{eqnarray*}
g_i(z,x) & = & g^{(i)}_{j_1j_2} dz^{j_1} \otimes dz^{j_2} +
g^{(i)}_{\al \be}dx^{\al} \otimes dx^{\be} + g^{(i)}_{j \al }dz^{j}
\otimes dx^{\al}
\end{eqnarray*}
and it is well known that in this coordinate system
\begin{eqnarray*}
g^{(i)}_{\al \be} = \delta_{\al \be}+ \bigo{|x|^2} & \qquad
\mbox{and} \qquad & g^{(i)}_{j \al } = \bigo{|x|}
\end{eqnarray*}\\

In order to perform the identification between $W_1$ and $W_2$ and
in order to glue the metrics together and define $g_\ep$, we
partially change the coordinate system, by setting \[ x = \ep \,
\ex^{-t} \, \theta \] on $F_1^{-1}(W_1)$ and \[ x = \ep \, \ex^{t}
\theta \] on $F_2^{-1}(W_2)$, for $\ep \in (0,1)$, $\log\ep < t < -
\log\ep$, $\theta \in S^{m-k-1}$. \\


Using these changes of coordinates the expressions of the two
metrics $g_1$ and $g_2$ on $U^k \times A^{1}_{\ep^2}$, where $
A^{1}_{\ep^2}$ is the annulus $\{\ep^2 < |x| < 1 \}$ become
respectively
\begin{eqnarray*}
g_{1}(z,t,\te) & = & g^{(1)}_{ij} dz^i \otimes dz^j \\
& + & u_{\ep}^{(1)} \, ^{\frac{4}{n-2}} \left[ \left(dt \otimes dt +
g^{(1)}_{\lambda \mu} d\theta^\lambda \otimes d\theta^\mu \right) +
g^{(1)}_{t \theta}dt \ltimes d\theta
\right]\\
& + & g^{(1)}_{i t} d\zi \otimes dt + g^{(1)}_{i \lambda} d\zi
\otimes d\te^{\lambda}
\end{eqnarray*}
and
\begin{eqnarray*}
g_{2}(z,t,\te) & = & g^{(2)}_{ij} dz^i \otimes dz^j \\
& + & u_{\ep}^{(2)} \, ^{\frac{4}{n-2}} \left[ \left(dt \otimes dt +
g^{(2)}_{\lambda \mu} d\theta^\lambda \otimes d\theta^\mu \right) +
g^{(2)}_{t \theta}dt \ltimes d\theta
\right]\\
& + & g^{(2)}_{i t} d\zi \otimes dt + g^{(2)}_{i \lambda} d\zi
\otimes d\te^{\lambda}
\end{eqnarray*}
where by the compact notation $g_{t \theta} \, dt \ltimes d\theta$
we indicate the general component of the normal metric tensor (that
is, it involves $dt \otimes dt$, $d\theta^\lambda \otimes
d\theta^\mu$
and $dt \otimes d\theta^\lambda$ components).\\

Remark that for $j = 1,2$ we have
\[
\begin{array}{rlllllll}
g^{(j)}_{\lambda \mu} & = & \bigo{1} & \qquad \qquad &
g^{(j)}_{t \theta} & =& \bigo{|x|^2}\\[3mm]
g^{(j)}_{i t} & = & \bigo{|x|^2} & \qquad \qquad &
g^{(j)}_{i \lambda} & = & \bigo{|x|^2}\\
\end{array}
\]
and
\begin{eqnarray*}
u_{\ep}^{(1)} (t) = \ep^{\frac{n-2}{2}} \ex^{-\frac{n-2}{2}t} &
\qquad \mbox{and} \qquad & u_{\ep}^{(2)} (t) = \ep^{\frac{n-2}{2}}
\ex^{\frac{n-2}{2}t}
\end{eqnarray*}

We choose a cut-off function $\chi : (\log\ep, -\log\ep) \rightarrow
[0,1]$ to be a non increasing smooth function which is identically
equal to $1$ in $(\log\ep, -1]$ and $0$ in $[1,-\log\ep)$ and we
choose another cut-off function $\eta : (\log\ep, -\log\ep)
\rightarrow [0,1]$ to be a non increasing smooth function which is
identically equal to $1$ in $(\log\ep, -\log\ep -1]$ and which
satisfies $\lim_{t\rightarrow -\log\ep} \eta = 0$. Using these two
cut-off functions, we can define a new normal conformal factor
$u_{\ep}$ by
\begin{eqnarray*}
u_{\ep} (t) : = \eta(t) \, u_{\ep}^{(1)} (t) + \eta(-t) \,
u_{\ep}^{(2)}(t)
\end{eqnarray*}
and the metric $g_{\ep}
$ by
\begin{eqnarray}
g_{\ep}(z,t,\te) & := & \left(\chi g^{(1)}_{ij} + (1-\chi) g^{(2)}_{ij} \right)dz^i \otimes dz^j \nonumber \\
& + & u_{\ep}^{\frac{4}{n-2}} \left[ dt \otimes dt + \left(\chi
g^{(1)}_{\lambda \mu} + (1-\chi) g^{(2)}_{\lambda \mu}
\right)d\theta^\lambda \otimes d\theta^\mu \right. \nonumber \\
& & \qquad \qquad \qquad + \left. \left(\chi g^{(1)}_{t \theta}
+ (1-\chi) g^{(2)}_{t \theta} \right) dt \ltimes d\theta \right]\\
& + & \left( \chi g^{(1)}_{i t} + (1-\chi) g^{(2)}_{i t} \right)
d\zi
\otimes dt \nonumber \\
& + & \left(\chi g^{(1)}_{i \lambda} + (1-\chi) g^{(2)}_{i
\lambda}\right) d\zi \otimes d\te^{\lambda} \nonumber\\ \nonumber
\end{eqnarray}

Closer inspection of this expression shows that the only objects
that are not {\it a priori} globally defined on the identification
of the tubular neighborhoods of $\iota_{1} (K)$ in $M_1$ and
$\iota_2 (K)$ in $M_2$ are the functions $\chi$ and $u_{\ep}$ (since
$\eta$ is used in the construction). However, observe that both
cut-off functions can easily be expressed as functions of the
Riemannian distance to $K$ in the respective manifolds. Hence they
are globally defined and the metric $g_{\ep}$ - whose definition can
be obviously completed by putting $g_{\ep} \equiv g_1$ and $g_{\ep}
\equiv g_2$ out of the "polyneck" - is a Riemannian metric which is
globally defined on the manifold $M$.

\section{Estimate of the scalar curvature}

Now we want to estimate the difference $S_{g_{\ep}} - S$ on the
"polyneck" (which, in the above coordinates, corresponds to $\log
\ep + 1 \leq t \leq -\log \ep -1$). To begin with, we restrict our
attention to the case where $\log \ep \leq t \leq -1$. Here the
normal conformal factor can be written down as $u_{\ep} =
u_{\ep}^{(1)} \left( 1 + {u_{\ep}^{(2)}}/{u_{\ep}^{(1)}} \right) $
so, if we define $h = u_{\ep}^{(2)} / u_{\ep}^{(1)}$ the metric
$g_{\ep}$ looks like
\begin{eqnarray*}
g_{\ep}(z,t,\te) = g^{(1)}_{ij} dz^i \otimes dz^j + \left( 1 + h
\right)^{\frac{4}{n-2}} g^{(1)}_{\al \be}dx^{\al} \otimes dx^{\be} +
g^{(1)}_{i \al} dz^{i} \otimes dx^{\al}
\end{eqnarray*}
where in fact $h = \ex^{(n-2)t} =
\ep^{(n-2)}|x|^{2-n}$.\\

In order to simplify the notations, let us drop the upper $^{(1)}$
indices and simply write
\begin{eqnarray*}
g(z,x,h) & = & g_{ij} dz^i \otimes dz^j + \left( 1 + h
\right)^{\frac{4}{n-2}} g_{\al \be}dx^{\al} \otimes dx^{\be} + g_{i
\al} dz^{i} \otimes dx^{\al}
\end{eqnarray*}
Recall that the following expansions hold
\begin{eqnarray*}
g_{ij} & = & g^K_{ij} (z) + \bigo{|x|} \\
g_{\al \be} & = & \delta_{\al \be}+ \bigo{|x|^2} \\
g_{i \al } & = & \bigo{|x|}
\end{eqnarray*}
In the following computation we will use the notations
\begin{eqnarray*}
g_h(z,x) & := & g(z,x,h) \\
g_0(z,x) & := & g(z,x,0) \\
\gtil_h(z) & := & g(z,0,h) \\
\gtil_0(z) & := & g(z,0,0) \\
\end{eqnarray*}
and their respective scalar curvature will be denoted by
\begin{eqnarray*}
S_h & := & S_{g_h}\\
S_0 & := & S_{g_0}\\
\tilde{S}_h & := & S_{\gtil_h}\\
\tilde{S}_0 & := & S_{\gtil_0}\\
\end{eqnarray*}

The idea is to estimate the difference between the scalar curvatures
of the metrics $g_h$ and $g_0$ by first estimating the differences
with the scalar curvature of the Riemannian product metrics
$\gtil_h$ and $\gtil_0$. In fact, we can easily obtain
\begin{eqnarray*}
\tilde{S}_h & = & \tilde{S}_0 + \left( 1 +
h\right)^{\frac{4}{n-2}}\Delta^{(x)}_{eucl}h
\end{eqnarray*}

Next we consider the term $S_h - \tilde{S}_h$. To keep notations
short, we agree that $A^{(j)}_k = A^{(j)}_l (z,x,h)$, $j,l \in
\mathbb{N}$ is a function, a row vector or a matrix whose
coefficients satisfy
\begin{eqnarray*}
\left| A^{(j)}_l (z,x,h) \right| & \leq & C \, |x|^l \\
\left| A^{(j)}_l (z,x,h) - A^{(j)}_k (z,x,0) \right| & \leq & C \, |x|^l \, |h| \\
\end{eqnarray*}
for some positive constant $C = C(j)$.\\

We start with the expansions of the coefficients of the metrics
$g_h$ (and hence also $g_0$ which corresponds to $g_h$ when $h=0$)
and their inverses in terms of $|x|$
\begin{eqnarray*}
g^{(h)}_{ij} & = & \gtil^{(h)}_{ij} (z) + \bigo{|x|} \\
g^{(h)}_{\al \be} & = & \gtil^{(h)}_{\al \be}+ \bigo{|x|^2} \\
g^{(h)}_{i \al } & = & \bigo{|x|}
\end{eqnarray*}
and
\begin{eqnarray*}
g^{ij}_{(h)} & = & \gtil^{ij}_{(h)} (z) + A_1^{(1)} \\
g^{\al \be}_{(h)} & = & \gtil^{\al \be}_{(h)} + A_1^{(2)} \\
g^{i \al }_{(h)} & = & A_1^{(3)}
\end{eqnarray*}

We estimate the Christoffel symbols of the metric $g_h$. Observe
that
\begin{eqnarray*}
g^{\ldots}_{(h)} \dd{g^{(h)}_{\ldots}}{\ldots} & = &
\left(\gtil^{\ldots}_{(h)} + A^{(4)}_1 \right) \left(
\dd{\gtil^{(h)}_{\ldots}}{\ldots} + \A{5} + \A{6}_1 \left[\nabla h
\right]\right) \\
& = & \gtil^{\ldots}_{(h)} \dd{\gtil^{(h)}_{\ldots}}{\ldots} +
\A{7}_0 + \A{8}_1 \left[ \nabla h \right]
\end{eqnarray*}
As a consequence we have that
\begin{eqnarray*}
\Gamma (h,\nabla h) & = & \tilde{\Gamma} (h, \nabla h) + \A{9}_0 +
\A{10}_1 \left[ \nabla h \right]
\end{eqnarray*}
Moreover, it is straightforward to check that
\begin{eqnarray*}
\tilde{\Gamma} (h, \nabla h) & = & \A{11}_0 + \A{12}_0 \left[ \nabla
h \right]
\end{eqnarray*}
Proceeding with the computation we get
\begin{eqnarray*}
\dd{\Gamma}{\ldots} (h,\nabla h) & = &
\dd{\tilde{\Gamma}}{\ldots}(h,\nabla h) + \A{13}_0 \left[ \nabla
h\right] + \A{14}_1 \left[ \nabla h , \nabla h\right] + \A{15}_1
\left[ \nabla^2 h\right]
\\
\\
\dd{\tilde{\Gamma}}{\ldots}(h,\nabla h)& = & \A{16}_0 \left[ \nabla
h\right]+ \A{17}_0 \left[ \nabla h , \nabla h\right] + \A{18}_0
\left[ \nabla^2 h\right]
\end{eqnarray*}
while for the product of Christoffel symbols, we get
\begin{eqnarray*}
\Gamma \, \Gamma (h,\nabla h) & = & \tilde{\Gamma} \, \tilde{\Gamma}
(h, \nabla h) + \A{19}_0 + \A{20}_0[\nabla h] + \A{21}_1 \left[
\nabla h , \nabla h\right]
\end{eqnarray*}
and hence we get for the coefficients of the curvature tensors
\begin{eqnarray*}
R(h,\nabla h, \nabla^2 h) & = & \tilde{R}(h,\nabla h, \nabla^2 h) +
\A{22}_0 + \A{23}_0 \left[ \nabla
h\right] \\
& + & \A{24}_1 \left[ \nabla h , \nabla h\right] + \A{25}_1 \left[
\nabla^2 h\right]
\\
\\
\tilde{R}(h,\nabla h, \nabla^2 h)& = & \A{26}_0 + \A{27}_0
\left[\nabla h\right] + \A{28}_0 \left[\nabla h,\nabla h\right] +
\A{29}_0 \left[\nabla^2 h\right]
\end{eqnarray*}
Finally, observing that
\begin{eqnarray*}
g^{\ldots}_h g^{\ldots}_h & = & \gtil^{\ldots}_h \gtil^{\ldots}_h +
\A{30}_1
\end{eqnarray*}
and contracting twice the Riemann tensor, we get the expression for
the scalar curvature
\begin{eqnarray*}
S_h & = & \tilde{S}_h + \A{31}_0 + \A{32}_0 \left[\nabla h\right] +
\A{33}_1 \left[ \nabla h , \nabla h\right] + \A{34}_1 \left[
\nabla^2 h\right]
\end{eqnarray*}

Choosing $h \equiv 0$ in the previous computation we obtain
immediately
\begin{eqnarray*}
S_0 & = & \tilde{S}_0 + \A{35}_0 (z,x,0)
\end{eqnarray*}
Hence we have obtained
\begin{eqnarray*}
S_h \quad = \quad S_0 & + & (1+h)^{-\frac{n+2}{n-2}}
\Delta^{(x)}_{eucl} h + \A{36}_0(z,x,h) - \A{36}_0 (z,x,0)
\\[3mm]
& + & \A{37}_0 \left[\nabla h\right] + \A{38}_1 \left[ \nabla h , \nabla h\right] + \A{39}_1 \left[ \nabla^2 h\right] \\
\end{eqnarray*}
Since $h = \ep^{n-2} |x|^{2-n}$ is $\Delta^{(x)}_{eucl}$-harmonic we
conclude that
\begin{eqnarray*}
S_h - S_0 & = & \A{40}_0 + \A{41}_0 [\nabla h ] + \A{42}_1 [\nabla h
, \nabla h] + \A{43}_1 [\nabla^2 h] \\
& = & \bigo{\ep^{n-2} |x|^{1-n}}\\
& = & \bigo{\ep^{-1} \ex^{(n-1)t}}
\end{eqnarray*}
We remark that, when $t = \log \ep + 1$, we get
the estimate $S_{g_{\ep}} - S_{g_1} = {\mathcal O} (\ep^{n-2})$. \\

Let us now treat the case where $-1 \leq t \leq 0$. The action of
the cut-off function is effective here, so {\it a priori} we have to
handle the full expression of $g_{\ep}$. In any case, it is easy to
see that one can always write for $-1 \leq t \leq 0$
\begin{eqnarray*}
g_{\ep}(z,t,\theta) & = & \left(g^{1}_{ij} + \bigo{|x|}\right) dz^i \otimes dz^j \\
& + & \left( 1 + h \right)^{\frac{4}{n-2}} \left( g^{(1)}_{\al \be}
+ \bigo{|x|}\right) dx^{\al} \otimes dx^{\be}\\
& + & \left( g^{(1)}_{i \al} + \bigo{|x|} \right) dz^{i} \otimes
dx^{\al}
\end{eqnarray*}
Hence, if we take $g(z,x,h) = g_{\ep}$ and $g(z,x,0) = g_{1} +
\bigo{|x|}$ in the previous computation we get immediately
$S_{g_{\ep}} - S_{g_1 + \bigo{|x|}} = \bigo{\ep^{n-2} |x|^{1-n}}$.\\

Now we observe that in general if we have two metrics $g$ and
$\hat{g}$ such that $\hat{g} = g + \bigo{|x|}$, then $\hat{\Gamma} =
\Gamma + \bigo{1}$ and $\hat{R} = R + \bigo{|x|^{-1}}$, so the
scalar curvatures of $g$ and $\hat{g}$ are related by $\hat{S}
= S + \bigo{|x|^{-1}}$.\\

To conclude, we have that $$ S_{g_{\ep}} - S_{g_1} = \bigo{|x|^{-1}}
= \bigo{\ep^{-1}\ex^{t}} $$ for $-1 \leq t \leq 0$. In particular,
when $t=0$ we get $S_{g_{\ep}} - S_{g_1} = \bigo{\ep^{-1}}$. Similar
estimates hold for $S_{g_{\ep}} - S_{g_{2}}$ when $0 \leq t \leq
-\log\ep -1$ and hence we have obtained the
\begin{lemma}
There exists a constant $c >0$ independent of $\ep \in (0,1)$ such
that
\[
|S_{g_{\ep}} - S | \leq c \, \ep^{-1} \, (\cip \, t)^{1-n}
\]
for $|t| \leq |\log \ep | -1$.
\end{lemma}

\section{Analysis of a linear operator}

In order to obtain the proof of the main Theorem, we want to solve,
using a perturbation argument, the Yamabe equation
\begin{eqnarray}
\label{Yamabe} \Delta_{\gep}u + c_n S_{\gep} u & = & c_n S
u^{\frac{n+2}{n-2}}
\end{eqnarray}
where $c_n = - (n-2) / 4 (n-1)$.\\

If we are able to find such a function $u$, then, by performing the
conformal change $\tilde{g}_{\ep} = u^{\frac{4}{n-2}}\gep$ we get a
metric $\tilde{g}_{\ep}$, whose scalar curvature is the constant equal to $S$.\\

We write $u = 1 + v$ where $v$ is a small function ($|v|\leq 1/2$)
so that the equation becomes
\begin{eqnarray}
\label{Yamabe in 1} \Delta_{\gep} v - \frac{4c_n}{n-2} S_{\gep} v &
= & c_n \left( S- S_{\gep} \right) + c_n \frac{n+2}{n-2}\left( S- S_{\gep} \right) v \\
& + & c_n S \left((1 + v)^{\frac{n+2}{n-2}} -1 -\frac{n+2}{n-2}v
\right)\nonumber
\end{eqnarray}

We define the linearized scalar curvature operator by
\[
L_{\gep} : =\Delta_{g_{\ep}} - \frac{4c_n}{n-2} S_{\gep}
=\Delta_{\gep} + \frac{S_{\gep}}{n-1}
\]
Our aim is to study the operator $L_{\gep}$ and provide an {\it a
priori} estimate for the solutions of the linear problem \[ L_{\gep}
v = f
\]
This is the starting point and the key-tool for the nonlinear perturbation argument.\\

Unfortunately a global {\it a priori} estimate is not immediately
available. We will be able to obtain such an estimate using an
argument by contradiction, once a local {\it a priori} estimate is
obtained for the solutions of the linearized problem on the "polyneck".\\

\subsection{Local expression for $\Delta_{\gep}$ on the "polyneck" and
barrier functions}

The first step is to write down the local expression for the
$\gep$-laplacian, which is the principal part of our operator, on
the "polyneck". Clearly, we can restrict ourselves to the set
$\{\log \ep + 1 \leq t \leq 0\}$ where $|x| =\ep \ex^{-t}$. We have
at hand the expansions
\begin{eqnarray*}
\geps_{ij} & = & g^K_{ij}(z) + \bigo{|x|}\\
\geps_{it} & = & \bigo{|x|^2}\\
\geps_{i \lambda} & = & \bigo{|x|^2}\\
\geps_{tt} & = & \uep^{\frac{4}{n-2}} \left(1 + \bigo{|x|^2}
\right)\\
\geps_{t\lambda} & = & \uep^{\frac{4}{n-2}} \bigo{|x|^2} \\
\geps_{\lambda\mu} & =& \uep^{\frac{4}{n-2}} \left( g^{}_{\lambda
\mu}(\theta) + \bigo{|x|^2} \right)
\end{eqnarray*}
where $g_{\lambda \mu} (\theta)$ is the common value of
$g^{(1)}_{\lambda \mu} (\theta)$ and $g^{(2)}_{\lambda \mu}
(\theta)$. Hence
\begin{eqnarray*}
\sqrt{\gep} & = & \sqrt{{\rm det}(g_{ij}^K(z))} \,\, \sqrt{{\rm
det}(g_{\lambda \mu}(z))} \,\, \uep^{\frac{2n}{n-2}}(t) \,\, \left[1
+ \bigo{|x|}\right]
\end{eqnarray*}
So, for coefficients of the inverse matrix we have the expansions
\begin{eqnarray*}
\gep^{ij} & = & g_K^{ij}(z) + \bigo{|x|}\\
\gep^{it} & = & \bigo{|x|^2}\\
\gep^{i \lambda} & = & \bigo{|x|^2}\\
\gep^{tt} & = & \uep^{-\frac{4}{n-2}} \left[ 1 + \bigo{|x|}\right]\\
\gep^{t\lambda} & = & \bigo{|x|^2} \\
\gep^{\lambda\mu} & = & \uep^{-\frac{4}{n-2}}g_{}^{\lambda \mu}
\left[ 1 + \bigo{|x|} \right]
\end{eqnarray*}
A straightforward computation yields the expression we were looking
for
\begin{multline}
\label{locallapl} \Delta_{\gep} = \uep^{-\frac{4}{n-2}} \left[
\partial_t^2 + (n-2) \, \tip \left(\frac{n-2}{2}t \right)
\partial_t +  \Delta^{(\theta)}_{S^{n-1}} + \uep^{\frac{4}{n-2}}\Delta^{(z)}_{K}
 + \bigo{|x|}\Phi(\nabla, \nabla^2)\right]
\end{multline}
where $\Phi(\nabla, \nabla^2)$ is a nonlinear differential operator
involving first order and second order partial derivatives with
respect to $t$, $\theta^\lambda$ and $z^j$ and whose coefficients
are bounded uniformly on the "polyneck",
as $\ep \in (0,1)$.\\

To obtain the local {\it a priori} estimates, the key tool is the
maximum principle for the ${\gep}$-Laplacian and the construction of
barrier functions. In order to find the later, let us remark that
\begin{eqnarray*}
\left( \partial_t^2 + \left( \frac{n-2}{2} \right)^2 \right) \left(
\cip \left( \frac{n-2}{2}t\right) u \right) & = & \left(\cip \left(
\frac{n-2}{2}t\right) \partial_t^2 + (n-2)\sip \left(
\frac{n-2}{2}t\right) \right) u
\end{eqnarray*}
So we can conjugate the $\gep$-Laplacian by a multiple of the
function $\cip(t (n-2) /2)$ - in particular, of course, by $\uep$ -
to obtain the following identity
\begin{eqnarray}
\label{Conjugation} \Delta_{\gep} & = & \uep^{-\frac{n+2}{n-2}} \,
\mathcal{L}_{\ep} \,\left( \uep \,\cdot \right)
\end{eqnarray}
where
\[
\lep \quad = \quad
\partial_t^2 \qu - \qu \left(\frac{n-2}{2}\right)^2 \qu
+ \qu \Delta^{(\theta)}_{S^{n-1}} \qu + \qu
\uep^{\frac{4}{n-2}}\Delta^{(z)}_{K} \qu + \qu \bigo{|x|}
\tilde\Phi(\nabla, \nabla^2)
\]
where the linear second order differential operator $\tilde
\Phi(\nabla, \nabla^2)$ enjoys similar properties as the operator
$\Phi$ above. For $(n-2)/2 \leq \delta \leq 0$ we have that
\begin{eqnarray*}
\lep (\cip t)^{\delta} & = & \left[ \delta^2 - \left( \frac{n-2}{2}
\right)^2 + \bigo{|x|} \right] (\cip t)^{\delta} + \left( \delta -
\delta^2 \right)
(\cip t)^{\delta -2}\\
\end{eqnarray*}
By our choice of the parameter $\delta $ we have immediately
\begin{eqnarray*}
\delta - \delta^2 \leq 0 \qquad \mbox{and} \qquad
\delta^2 - \left( \frac{n-2}{2} \right)^2 \leq 0 \\
\end{eqnarray*}
In order to estimate the term $\bigo{|x|}$ let us take $\al
> 0$ and let $\ep_{\al} \in (0,1)$ be chosen so that $\log \ep_{\al}
+ \al <0$ or equivalently $\ep_{\al} \ex^{\al} < 1$, then it is easy
to see that $|x| \leq \ex^{-\al}$ for every $\ep \in (0,\ep_{\al})$
and every $t \in [\log \ep + \al , 0 ]$. Finally, by choosing $\al >
0$ such that
\begin{eqnarray*}
\ex^{-\al} & \leq & - \frac{1}{2} \left( \delta^2 - \left(
\frac{n-2}{2} \right)^2 \right)
\end{eqnarray*}
we obtain that, for every $\ep \in (0, \ep_{\al})$ and for $t \in
[\log \ep + \al , 0 ]$
\begin{eqnarray*}
\lep (\cip t)^{\delta} & \leq & \frac{1}{2} \left( \delta^2 - \left(
\frac{n-2}{2} \right)^2 \right) (\cip t)^{\delta}
\end{eqnarray*}
When $0 \leq \delta \leq (n-2)/2$ we use the function $\cip ( \delta
t)$ and we get
\begin{eqnarray*}
\lep \cip (\delta t) & = & \left( \delta^2 - \left( \frac{n-2}{2}
\right)^2 + \bigo{|x|} \right) \cip \delta t \\
& \leq & \frac{1}{2} \left( \delta^2 - \left( \frac{n-2}{2}
\right)^2 \right) \cip \delta t
\end{eqnarray*}
with similar restrictions on $\ep$ and $t$.\\

We define the function $\varphi_{\delta}$ by
\begin{eqnarray*}
\varphi_{\delta}\quad = \quad \uep^{-1}(\cip t)^{\delta} & & \qquad
{\rm if} \quad \frac{n-2}{ 2} \leq \delta \leq 0
\\
\varphi_{\delta} \quad = \quad \uep^{-1}\cip \delta t \,\,\,\,& &
\qquad {\rm if} \quad 0 \leq \delta \leq \frac{n-2}{2}
\end{eqnarray*}
and taking into account the conjugation (\ref{Conjugation})
described above, we can state the following
\begin{lemma}
\label{MAX PRINC} Given $ \delta \in (-\frac{n-2}{ 2},
\frac{n-2}{2})$ there exist a real number $\al = \al (n, \delta) >
0$ and a constant $C = C(n, \delta) \geq 0$ such that for every $\ep
\in (0, \ep_{\al})$ we have
\begin{eqnarray}
\label{max princ} \Delta_{\gep} \varphi_{\delta} & \leq & - C
\uep^{-\frac{4}{n-2}} \varphi_{\delta}
\end{eqnarray}
in the set $T^{\ep}_{\al} = \{\log \ep + \al \leq t \leq -\log \ep -
\al \}$. \label{le:4.1}
\end{lemma}
In particular the functions $\varphi_{\delta}$ can be used as
barrier functions in the set $T^{\ep}_{\al} = \{\log \ep + \al \leq
t \leq -\log \ep - \al \}$.

\subsection{Local {\it a priori} estimate using the maximum principle}

We first provide a local {\it a priori} estimate for the
$g_{\ep}$-Laplacian, then we will observe that a similar estimate
holds for the operator $L_{\gep}$. This later estimate uses the
scalar curvature estimate of the previous section since the term
$S_{\gep}$
appears in the expression of $L_{\gep}$.\\

Let us assume that $v , f$ are bounded functions satisfying
$\Delta_{\gep} v = f $ in $T^{\ep}_{\al}$. The inequality found in
Lemma~\ref{MAX PRINC} multiplied by a nonnegative real constant $a
\geq 0$ yields
\begin{eqnarray*}
\Delta_{\gep} \left( a \varphi_{\delta} - v \right) & \leq & -a C
\uep^{-\frac{4}{n-2}} \varphi_{\delta} - f
\end{eqnarray*}
If we chose
\begin{eqnarray*}
a = C' \left(\,\, \sup_{T^{\ep}_{\al}} \left|\uep^{\frac{4}{n-2}}
\varphi_{\delta}^{-1} f \right|\,\, + \,\, \sup_{\partial
T^{\ep}_{\al}} \left| \varphi_{\delta}^{-1} v \right| \,\,\right)
\end{eqnarray*}
where $C' = \max \{ 1,C^{-1} \}$ and $\partial T^{\ep}_{\al} = \{ t
= \pm \log \ep \pm \al \} $, we obtain immediately
\begin{eqnarray*}
\Delta_{\gep} \left( a \varphi_{\delta} - v \right) & \leq & 0
\qquad {\rm in} \quad T^{\ep}_{\al} \\
a \varphi_{\delta} - v & \geq & 0 \qquad {\rm on} \quad \partial
T^{\ep}_{\al}
\end{eqnarray*}
Hence, by the maximum principle $a \varphi_{\delta} - v \geq 0$ on
 $T^{\ep}_{\al}$. In particular, we get
\begin{eqnarray*}
\sup_{T^{\ep}_{\al}} \left| \varphi_{\delta}^{-1} v \right| & \leq &
C' \left(\,\, \sup_{T^{\ep}_{\al}} \left|\uep^{\frac{4}{n-2}}
\varphi_{\delta}^{-1} f \right|\,\, + \,\, \sup_{\partial
T^{\ep}_{\al}} \left| \varphi_{\delta}^{-1} v \right| \,\,\right)
\end{eqnarray*}
In order to simplify the above expression, which is the estimate we
were looking for, it is sufficient to replace $\uep$ by its
expression and to observe that for every $\lambda \in \mathbb{R}$
there exist two constants $K_1({\lambda}),K_2({\lambda}) \geq 0$
such that
$$K_1({\lambda}) \left( \cip t \right)^{\lambda} \leq \cip \lambda t
\leq K_2({\lambda}) \left( \cip t \right)^{\lambda}$$ for $t
\in {\mathbb R}$. \\

Performing simple manipulations, the above estimate can be written
as
\begin{eqnarray*}
\sup_{T^{\ep}_{\al}} \left| \psi_{\ep}^{\frac{n-2}{2}-\delta} v
\right| & \leq & C_{n, \delta} \left(\,\, \sup_{T^{\ep}_{\al}}
\left| \psi_{\ep}^{\frac{n+2}{2}-\delta} f \right|\,\, + \,\,
\sup_{\partial T^{\ep}_{\al}} \left|
\psi_{\ep}^{\frac{n-2}{2}-\delta} v \right| \,\,\right)
\end{eqnarray*}
where $\psi_{\ep} = \ep \, \cip t$.\\

Now let us assume that $v,f \in {\mathcal
C}^{\infty}(T^{\ep}_{\al})$ are functions verifying $L_{\gep} v =
f$. By the previous result we immediately have
\begin{eqnarray*}
\sup_{T^{\ep}_{\al}} \left| \psi_{\ep}^{\frac{n-2}{2}-\delta} v
\right| & \leq & C_{n, \delta} \left(\,\, \sup_{T^{\ep}_{\al}}
\left| \psi_{\ep}^{\frac{n+2}{2}-\delta} f \right|\,\, + \,\,
\sup_{T^{\ep}_{\al}} \left| \psi_{\ep}^{\frac{n+2}{2}-\delta}\,
 S_{\gep} v \right|\,\, + \,\, \sup_{\partial T^{\ep}_{\al}} \left|
\psi_{\ep}^{\frac{n-2}{2}-\delta} v \right| \,\,\right)
\end{eqnarray*}
Now let us look at the term $\left|
\psi_{\ep}^{\frac{n+2}{2}-\delta} \, S_{\gep} v \right| $, which can
be obviously written as $\left|\psi_{\ep}^2 S_{\gep}\right|
\left|\psi_{\ep}^{\frac{n-2}{2}-\delta}v\right| $. The only term to
control is the factor $\left|\psi_{\ep}^2 S_{\gep}\right|$, but
thanks to the scalar curvature estimate we can say that, for a
suitable constant $C''>0$
\begin{eqnarray*}
\left|\psi_{\ep}^2 S_{\gep}\right| \leq C'' \, (\ep + \ex^{-\al})
\end{eqnarray*}
for all $\ep \in (0, \ep_{\al})$. Hence, if $\al>0$ is fixed large
enough, we get
\begin{eqnarray*}
C_{n, \delta} \, \sup_{T^{\ep}_{\al}} \left|
\psi_{\ep}^{\frac{n+2}{2}-\delta} \, S_{\gep} v \right| & \leq &
\frac{1}{2}\sup_{T^{\ep}_{\al}}\left|\psi_{\ep}^{\frac{n-2}{2}-\delta}
v \right|
\end{eqnarray*}
Introducing this information back in the above estimate, we get
\begin{prop}
\label{LAPE} Given $ \delta \in (-\frac{n-2}{ 2}, \frac{n-2}{2})$,
there exist a real number $\al = \al (n, \delta) > 0$ and a constant
$C_{n, \delta} \geq 0$ such that for all $\ep \in (0, \ep_{\al})$
and all $v,f \in {\mathcal C}^{0}(T^{\ep}_{\al})$ satisfying
$L_{\gep} v = f$, the following estimate holds
\begin{eqnarray}
\label{lape} \sup_{T^{\ep}_{\al}} \left|
\psi_{\ep}^{\frac{n-2}{2}-\delta} v \right| & \leq & C_{n, \delta}
\left(\,\, \sup_{T^{\ep}_{\al}} \left|
\psi_{\ep}^{\frac{n+2}{2}-\delta} f \right|\,\, + \,\,
\sup_{\partial T^{\ep}_{\al}} \left|
\psi_{\ep}^{\frac{n-2}{2}-\delta} v \right| \,\,\right)
\end{eqnarray}
where $\psi_{\ep} = \ep \, \cip t$.
\end{prop}

\subsection{Global {\it a priori} estimate}

Thanks to the previous local result, we will be able to prove a
global {\it a priori} estimate. To introduce the result, we define a
smooth function $\psi_{\ep}$ by
\[
\psi_{\ep} : =
\begin{cases}
\ep \cip t & \text{in $T^{\ep}_{\al}$} \\
1 & \text{in $M \setminus T^{\ep}_0$}
\end{cases}
\]
where $T^{\ep}_{\rho} : = \{\log \ep + \rho \leq t \leq - \log \ep
-\rho \}$, for $\rho > 0$ and $\psi_\ep$ interpolate smoothly
between these definitions in $T^{\ep}_0 \setminus T^{\ep}_{\al}$.

\begin{prop}
\label{GAPE} Let $M = M_1 \, \sharp_{K} \, M_2$ be the generalized
connected sum obtained by removing an $\ep$-tubular neighborhood
$V^{\ep}_i$ of $\iota_i(K)$ from each $M_i$, $i= 1,2$ and
identifying the two boundaries. Suppose that both $L_{g_1}$ and
$L_{g_2}$ have trivial kernels on $M_1$ and on $M_2$ respectively,
then for every $ \delta \in \left(-\frac{n-2}{ 2},
\frac{n-2}{2}\right)$ there exist a real number $\al = \al (n,
\delta) > 0$ and a constant $C_{n, \delta} \geq 0$ such that for
every $\ep \in (0, \ep_{\al})$ and every functions $v,f \in
{\mathcal C}^{0}(M)$ satisfying $L_{\gep} v = f$, the following
estimate holds
\begin{eqnarray}
\label{gape} \sup_{M} \left| \psi_{\ep}^{\frac{n-2}{2}-\delta} v
\right| & \leq & C_{n, \delta} \left( \,\, \sup_{M} \left|
\psi_{\ep}^{\frac{n+2}{2}-\delta} f \right| \right)
\end{eqnarray}
\end{prop}
The proof is by contradiction. Let us assume that the statement is
false. Then for every $j \in \mathbb{N}$ we can find a triple
$(\ep_j , v_j , f_j)$ such that
\begin{enumerate}
  \item $\epj < \ex^{-j}$
  \item $L_{g_{\ep_j}} v_j = f_j$
  \item $\sup_{M} \left| \psi_{\ep_j}^{\frac{n-2}{2} - \delta} v_j
  \right| = 1 $
  \item $ \lim_{j \rightarrow \infty} \sup_{M} \left| \psi_{\ep_j}^{\frac{n+2}{2}} f_j \right| = 0 $
\end{enumerate}

For every $j \in \mathbb{N}$ we consider a point $p_j$ such that
$\left| \psi_{\ep_j}^{\frac{n-2}{2} - \delta} (p_j) v_j(p_j) \right|
= 1 $, then (up to a subsequence) we have to distinguish two cases :

\begin{description}
  \item[Case 1]
  $p_j \in M \setminus T^{\ep_j}_{\al}$ for every $j\in \mathbb{N}$
  \item[Case 2]
  $p_j \in T^{\ep_j}_{\al}$ for every $j \in \mathbb{N}$
\end{description}

Without loss of generality we can assume (up to a subsequence) that
$p_j \in M_1 \setminus V^{\ep_j}_1$, for all $j \in \mathbb{N}$, so,
in the first case all the $p_j$'s are in the compact set $
Q_1^{\ex^{-\al}} = M_1 \setminus V^{\ex^{-\al}}_1$, then (up to a
subsequence) they must converge to a point $p_{\infty} \in
Q_1^{\ex^{-\al}}$. We prove now that, for every compact set
$Q^{\sigma} = Q_1^{\sigma} \cup Q_2^{\sigma} = \left(M_1 \setminus
V^{\sigma}_1 \right) \, \cup \, \left(M_2 \setminus
V^{\sigma}_2\right)$, $\sigma
>0$, the sequence of functions $\{ v_j \}_{j \in \mathbb{N}}$
converges (up to a subsequence) to a function $v_{\infty}$ in
$L^\infty(Q^{\sigma})$. This will in particular imply that
$|v_{\infty}
(p_{\infty}) | > 0$.\\

In order to prove the uniform convergence of the $v_j$'s on the
compact $Q^{\sigma}$, we start by observing that $$|v_j| \leq
\left(\inf_{Q^{\sigma}} \left|\psi_j^{\frac{n-2}{2} - \delta}
\right| \right)^{-1} \leq \frac{2}{\sigma}$$ and hence
$\nor{v_j}{L^{\infty}(Q^{\sigma})} \leq 2 / \sigma$. \\

The next step is to get a ${L^{\infty}(Q^{\sigma})}$- uniform bound
for $\nabla v_j$. To do that we need the following $L^p$-regularity
result \cite{Gil-Tru} for solutions of linear elliptic equations
\begin{teor}
Let be $L = a^{ij} \partial_{ij} + b^i \partial_i + c $, where the
$a,b,c$'s are functions defined on an open domain $\Omega
 \subset \mathbb{R}^m$, let be $1<p<\infty$ and let be $u \in W^{2,p}_{loc} (\Omega) \cap
 L^p(\Omega)$.
 Moreover suppose that:
\begin{enumerate}
  \item $a^{ij} \in \mathcal{C}^0 (\Omega)$; $b^j,c \in L^{\infty
  }(\Omega)$; $f \in L^p(\Omega)$
  \item There exist $\lambda, \Lambda >0 $ such that $|a^{ij}|,|b^j|,|c| \leq \Lambda$ and $a^{ij} \xi_i \xi_j \geq \lambda
  |\xi|^2$ for every $\xi \in \mathbb{R}^n$
  \item $Lu = f$
\end{enumerate}
then, for every $\Omega' \subset \subset \Omega$, the following
estimate holds:
\begin{eqnarray*}
\nor{u}{W^{2,p}(\Omega ')} & \leq & C \left( \nor{u}{L^p(\Omega)} +
\nor{f}{L^p(\Omega)} \right)
\end{eqnarray*}
for a suitable constant $C$.
\end{teor}
\noindent This result can be restated in our context by saying:
\begin{cor}
\label{GTcor} Let be $\sigma>0$ and suppose that the linear elliptic
differential operator $L_g = \Delta_g + c$ is defined on a geodesic
ball $B_{\sigma/2}$ of the Riemannian manifold $(M,g)$, where $c$ is
a continuous bounded function on $B_{\sigma/2}$. Moreover let be
$1<p<\infty$ and let be $u \in W^{2,p}_{loc} (B_{\sigma/2}) \cap
 L^p(B_{\sigma/2})$, $f \in L^p(B_{\sigma/2})$ such that $L_g u =
 f$, then for every $0<r< \sigma/2$ the following estimate holds
\begin{eqnarray*}
\nor{u}{W^{2,p}(B_r )} & \leq & C \left( \nor{u}{L^p(B_{\sigma/2})}
+ \nor{f}{L^p(B_{\sigma/2})} \right)
\end{eqnarray*}
for a suitable constant $C$ (depending on $\sigma$).
\end{cor}
\noindent In our case it is convenient to cover the compact set
$Q^{\sigma}$ by finitely many geodesic balls of radius $r = \sigma /
4$. We can state
\begin{eqnarray*}
\nor{v_{j}}{W^{2,p}(B_{\sigma/4})} & \leq & C \left(
\nor{v_{j}}{L^p(B_{\sigma/2})} + \nor{f_j}{L^p(B_{\sigma/2})} \right)\\
& \leq & C' \left(
\nor{v_{j}}{L^{\infty}(B_{\sigma/2})} + \nor{f_j}{L^{\infty}(B_{\sigma/2})} \right)\\
& \leq & \frac{C''}{\sigma}
\end{eqnarray*}
Thanks to Sobolev Embedding Theorem with $p>m/2$, $W^{2,p}$ is
continuously embedded in $L^{\infty}$, so $\nor{\nabla
v_j}{L^{\infty}(B_{\sigma/4})} \leq C'''/ \sigma$. Now, by Ascoli's
Theorem, we conclude that (up to a subsequence) the sequence $\{ v_j
\}_{j \in \mathbb{N}}$ converges uniformly to a function
$v_{\infty}$ on every $B_{\sigma/4}$. Using a classical diagonal
argument we have the convergence on each $Q^\sigma$.\\

To summarize, in the $\mathbf{Case \,\, 1}$, we have found a
subsequence such that $v_j \rightarrow v_{\infty}$ with respect to
the $L^{\infty}$-norm on any $Q^{\sigma}$, in particular $v_{\infty}
\in \mathcal{C}^0({Q^{\sigma}})$, and, for $\sigma = \ex^{-\al}$, we
get $|v_{\infty}(p_{\infty})| >0$, as we have already
remarked.\\

Now, let us consider $\mathbf{Case \,\, 2}$. Since each $p_j$ is in
$T^{\ep_j}_{\al}$, we can apply the local {\it a priori} estimate
(\ref{lape}) obtained in the previous section to get
\begin{eqnarray*}
C_{n, \delta}^{-1} & \leq & \sup_{T^{\ep_j}_{\al}} \left|
\psi_{\ep_j}^{\frac{n+2}{2}-\delta} f_j \right|\,\, + \,\,
\sup_{\partial T^{\ep_j}_{\al}} \left|
\psi_{\ep_j}^{\frac{n-2}{2}-\delta} v_j \right| \\
& \leq & \sup_{M} \left| \psi_{\ep_j}^{\frac{n+2}{2}-\delta} f_j
\right|\,\, + \,\, \max_{\partial V_1^{\ex^{-\al}}} \left|
\psi_{\ep_j}^{\frac{n-2}{2}-\delta} v_j \right| + \,\,
\max_{\partial V_2^{\ex^{-\al}}} \left|
\psi_{\ep_j}^{\frac{n-2}{2}-\delta} v_j \right|
\end{eqnarray*}
This shows that we can choose a sequence of points $q_j \in
\partial V_1^{\ex^{-\al}} \cup \partial V_2^{\ex^{-\al}} $ such that $$\lim_{j \rightarrow \infty}
\left| \psi_{\ep_j}^{\frac{n-2}{2}-\delta}(q_j) v_j(q_j) \right| =
C_{n, \delta}^{-1}$$

In particular we have that $\lim_{j \rightarrow\infty} \left|
v_{j}(q_j) \right| = 2 C_{n, \delta}^{-1}\ex^{\al}$, then by using
the $L^{\infty}$-convergence to $v_{\infty}$ on the compact set
$Q^{\ex^{-\al}}$, it is easy to see that $|v_{\infty} (q_{\infty})|
> 0$, where $q_{\infty} \in \partial
V_1^{\ex^{-\al}} \cup \partial V_2^{\ex^{-\al}}$ is the limit (up to
a subsequence) of the sequence $\{q_j\}$.\\

Hence, in both the cases, we have found a point $ P \in M \setminus
T^{\ep}_{\al}$ such that $v_{\infty} (P) \neq 0 $. Without loss of
generality we can suppose that $P \in M_1 \setminus \iota_1(K)$: if
we prove that $L_{g_1} v_{\infty} = 0 $ on $M_1$, then by the
hypothesis on the kernel of $L_{g_1}$, $v_{\infty}$ must be
identically zero and we have a contradiction.\\

Hence, it remains to prove that $v_{\infty }$ is in the kernel of
$L_{g_1}$. This will be achieved in two steps. The first one amounts
to say that $L_{g_1} v_{\infty} = 0 $ on $M_1 \setminus \iota_1(K)$
in the sense of distributions, the second one amounts to estimate
the growth of $v_{\infty}$ near $\iota_1(K)$ and then to conclude by
means of the following classical result.
\begin{prop}
Suppose that
\[
\left\{
  \begin{array}{rllllll}
    L_{g_1} \,u & = & 0 & \mbox{in} \, \mathcal{D'}(M_1
\setminus \iota_1(K)) \\[3mm]
    |u| & \leq & C \, \left|{\rm d}_{g_1} (\cdot , \iota_1(K) )\right|^{-\gamma} & \mbox{in} \, V^{\rho}_1
  \end{array}
\right.
\]
For $0< \gamma < n-2$, a suitable real number $\rho > 0$ and a
constant $C\geq 0$, then $u \in \mathcal{C}^{\infty} (M_1)$ and
satisfies $L_{g_1}u = 0$ on $M_1$.
\end{prop}
We choose $\varphi \in \mathcal{D}(M_1 \setminus \iota_1(K))$ and
$\sigma >0$ such that ${\rm supp} \, \varphi \subset Q^{\sigma}$. We
claim that
\begin{eqnarray*} \int_{M_1} \vin L_{g_1} \fhi
\dvol{g_1} & = & 0
\end{eqnarray*}
This identity is obtained by taking the limit, as $\epj$ tends to
$0$ in the expression
\begin{eqnarray*}
\int_{M} v_j \, L_{g_{\epj}}\fhi \dvol{g_{\epj}} = \int_{M} f_j \,
\fhi \dvol{g_{\epj}}
\end{eqnarray*}
Clearly, the right hand side of this expression tends to zero as
$\ep_j$ tends to $0$. As far as the right hand side is concerned
$g_{\epj}$ converges (in ${\mathcal C}^2$ topology) to $g_1$ on
$Q^\sigma$ and hence $L_{g_{\epj}} \fhi$ converges to $L_{g_1} \fhi$
in this set so that the left hand side converges to the required
expression as $\ep_j$ tends to $0$.\\

Finally we have to control the growth of $\vin$ near $\iota_1(K)$.
We remark that, on $V^{\rho}_1$
$$\frac{1}{2} |x| \leq \psi_{\ep_j} \leq 2 \, |x|$$
for every $j$, in particular
\begin{eqnarray*}
|x|^{\frac{n-2}{2} - \delta} |v_j| \leq C
\end{eqnarray*}
Hence $$|\vin| \leq C |x|^{\delta - \frac{n-2}{2}} = C
|x|^{-\gamma}$$ where $\gamma = \frac{n-2}{2} - \delta$. Since
$-\frac{n-2}{2} < \delta < \frac{n-2}{2}$, then $0 < \gamma < n-2$,
as needed.

\section{The nonlinear fixed point argument}

We are now ready to solve equation (\ref{Yamabe in 1}). Observe
that, as a consequence of the Proposition \ref{GAPE} , the operator
$L_{\gep}$ is injective for sufficiently small $\ep$. Since it is
also self-adjoint, then it is invertible. Now we are looking for a
function $v_{\ep} \in L^{\infty}(M)$, $\| v\|_{L^\infty (M)} \leq
1/2$ such that
\begin{eqnarray}
v_{\ep} = L_{\gep}^{-1} \circ F_\ep (v_{\ep})
\end{eqnarray}
where
$$
F_\ep(v) := c_n(S -S_{\gep}) + c_n(S -S_{\gep})v + c_n \left((1 +
v)^{\frac{n+2}{n-2}} -1 - \frac{n+2}{n-2} v \right)
$$
In other words we are looking for a fixed point for the operator
$L_{\gep}^{-1} \circ F_\ep$.\\

We claim that, for a suitable choice of $\delta$ and for
sufficiently small ${\ep}$ there exists a real number $r_{\ep} > 0$
such that
\begin{eqnarray*}
\nor{v}{L^{\infty}(M)} \leq r_{\ep} & \Longrightarrow &
\nor{L_{\gep}^{-1} \circ F_\ep (v)}{L^{\infty}(M)} \leq r_{\ep}
\end{eqnarray*}
Indeed, using the scalar curvature estimates it is easy to see that
\begin{eqnarray*}
\sup_{M} \left|\psi_{\ep}^{\frac{n+2}{2} - \delta} F_\ep (v) \right|
& \leq & C' \left( \ep^{n-2} + \ep^{\frac{n}{2} - \delta} +
\nor{v}{L^{\infty}(M)}^{2} \right)
\end{eqnarray*}
Now
\begin{eqnarray*}
\psi_{\ep}^{\delta - \frac{n-2}{2}} \left| \ep^{n-2} +
\ep^{\frac{n}{2} - \delta} + \nor{v}{L^{\infty}(M)}^{2} \right| &
\leq & C'' \left( \ep^{\frac{n-2}{2} + \delta} + \ep + \ep^{\delta -
\frac{n-2}{2}} \nor{v}{L^{\infty}(M)}^{2} \right)
\end{eqnarray*}
Therefore, using the estimate (\ref{gape}) and the hypothesis of the
claim we get
\begin{eqnarray*}
\nor{L_{\gep}^{-1} \circ F_\ep(v)}{L^{\infty}(M)} & \leq & C'''
\left( \ep^{\frac{n-2}{2} + \delta} + \ep + \ep^{\delta -
\frac{n-2}{2}} r_{\ep}^{2} \right)
\end{eqnarray*}
where $C''' = C C' C''$. To prove the claim it is sufficient to
choose $r_{\ep} >0$ such that
\begin{eqnarray*}
\ep^{\delta - \frac{n-2}{2}} r_{\ep}^{2} \leq r_{\ep}/(2C''') \qquad
\mbox{and} \qquad \ep^{\frac{n-2}{2} + \delta} + \ep \leq
r_{\ep}/(2C''')
\end{eqnarray*}
The first condition is satisfied if we choose $r_{\ep}=
\ep^{\frac{n-2}{2} - \delta} / (2C''') $. By this choice the second
inequality becomes
\begin{eqnarray*}
\ep^{2\delta} + \ep^{\delta - \left(\frac{n-2}{2} -1 \right)} \leq
1/(2C''')^2
\end{eqnarray*}
Now it is clear that if $\max\{0, (n-2)/ 2 \,-1 \} < \delta < (n-2)
/ 2 $, then it is possible to find $\ep_0 \in (0, \ep_{\al})$ such
that the last inequality is verified for all $\ep \in (0, \ep_0)$.
For those $\ep$'s, we can choose $r_{\ep} = \ep^{\frac{n-2}{2} -
\delta} / (2C''')$ and the claim follows, hence
\[
\nor{L_{\gep}^{-1} \circ F_\ep(v)}{L^{\infty}(M)} \leq r_{\ep}
\]

It is easy to check that the mapping \[ v \in L^\infty(M)
\longmapsto L_{\gep}^{-1} \circ F_\ep(v) \in L^\infty(M)
\]
is continuous and compact. This later property follows from the fact
that the equation we want to solve is a semilinear equation and
hence, if $v \in L^\infty (M)$ then $L^{-1}_{g_\ep} \circ F_\ep(v)
\in W^{2,p} (M)$ for all $p >1$. The claim follows from the fact
that the embedding $W^{2,p}(M) \longrightarrow L^{\infty} (M)$ is
compact, provided $p > m/2$. Applying Schauder's fixed point Theorem
yields the existence of a fixed point $v_\ep \in L^\infty (M)$ to
\[
v_\ep = L^{-1}_{g_\ep} \circ F_\ep(v_\ep)
\]
which satisfies $ \nor{v_\ep}{L^{\infty}(M)} \leq r_{\ep}$.\\

{\it A priori} the function $v_\ep$ is only bounded but, by a simple
boot-strap argument (based on Corollary \ref{GTcor}), one can easily
checks
that $v_\ep \in \mathcal{C}^{\infty} (M)$.\\

Finally, observe that as $\ep \rightarrow 0$, then $r_{\ep}
\rightarrow 0$ and consequently so does $\nor{v_\ep
}{L^{\infty}(M)}$. This shows that the conformal factor $u_\ep = 1 +
v_\ep$ is as close to $1$ as we want. This completes the proof of
the main Theorem. The estimate in the statement of the Theorem
follows at once from the definition of $r_\ep$.

\end{document}